\documentclass[12pt]{article}
\usepackage{amsmath,amsfonts,amsthm,amscd}
\newtheorem{theorem}{Theorem}

\newtheorem{proposition}{Proposition}
\newtheorem{lemma}{Lemma}
\newtheorem{definition}{Definition}
\newtheorem{corollary}{Corollary}

\numberwithin{equation}{section}
\numberwithin{theorem}{section}
\numberwithin{proposition}{section}
\numberwithin{lemma}{section}
\numberwithin{claim}{section}
\numberwithin{corollary}{section}
\begin{document}
\title{Fully nonlinear equations
on Riemannian manifolds with negative curvature} 
\author{ {Matthew J. Gursky}
\and 
{Jeff A. Viaclovsky}\thanks{Supported in
part by an NSF Postdoctoral Fellowship.}}
\date{May 1, 2002}
\maketitle
\section{Introduction}
Let $(M^n,g)$ be a compact, connected, $n$-dimensional  
Riemannian manifold, $n \geq 3$, and let the Ricci 
tensor and scalar curvature be denoted by  
$Ric$ and $R$, respectively. 
For a symmetric tensor $A$, we let $ \mbox{det}(A)$ denote 
the determinant of $A$, that is, the product 
of the eigenvalues of $A$. 
\begin{theorem}  
\label{det1}
Assume that $(M^n,g)$ has negative Ricci curvature. 
Then there exists a unique conformal metric $\tilde{g} 
= e^{2u} g$ with negative Ricci curvature 
satisfying
\begin{align}  
\mbox{\em{det}}(Ric_{\tilde{g}}) = \mbox{\em{constant}}.
\end{align}
\end{theorem}
By results of \cite{Brooks}, \cite{GaoYau}, and \cite{Lohkamp}, 
every compact manifold of dimension $n \geq 3$ admits a metric
with negative Ricci curvature. Therefore we have 
\begin{corollary}
Every smooth compact $n$-manifold, $n \geq 3$, admits a Riemannian metric with
$Ric < 0$ and  
\begin{align}
\mbox{\em{det}}(Ric) = \mbox{\em{constant}}. 
\end{align}
\end{corollary}
 This theorem may be viewed as a Monge-Amp\`ere version 
of a theorem of Aubin for the scalar curvature (\cite{Aubin}). 
Theorem \ref{det1} is a special case of a more general theorem
involving symmetric functions of eigenvalues of $Ric$, which we
describe next.
Let $A$ be a symmetric $n \times n$ matrix, and let 
$\sigma_k(A)$ denote the $k$th elementary symmetric 
of $A$. 
\begin{definition}
Let $(\lambda_1, \dots, \lambda_n) \in \mathbf{R}^n$.
We view the elementary symmetric functions as
functions on $\mathbf{R}^n$
$$\sigma_k(\lambda_1, \dots, \lambda_n) = \sum_{i_1 < \dots < i_k}
\lambda_{i_1} \cdots \lambda_{i_k},$$
and we define
$$\Gamma_k^+ = \mbox{component of } \{\sigma_k > 0\}
\mbox{ containing the positive
cone}.$$
We also define $\Gamma_k^- = - \Gamma_k^+$.
\end{definition}
For a symmetric linear transformation $A : V \rightarrow V$, where
$V$ is an $n$-dimensional inner product space, the
notation $A \in \Gamma_k^{\pm}$ will mean that
the eigenvalues of $A$ lie in the corresponding set.
We note that this notation also makes sense for a symmetric tensor on a
Riemannian manifold.
If $A \in \Gamma_k^+$, let
$\sigma_k^{1/k}(A) = \{ \sigma_k(A) \}^{1/k}$, 
and if $A \in \Gamma_k^-$, let
$\sigma_k^{1/k}(A) = - |\sigma_k(A)|^{1/k}$.

Let $t \in \mathbf{R}$, and define 
\begin{align}
A^t = \frac{1}{n-2} \left( Ric - \frac{t}{2(n-1)}Rg \right).
\end{align}
Note that for $t=1$, $A^1$ is the classical Schouten tensor (\cite{Eisenhart}).

The following is our main theorem.
\begin{theorem}
\label{sharpthm}
Assume that $A^t_g \in \Gamma_k^{-}$
for some $t<1$, and let $f(x) < 0$ be any smooth 
function on $M^n$. Then there 
exists a unique conformal metric $\tilde{g} = 
e^{2w}g$ satisfying
\begin{align}
\label{eqn1}
\sigma_k^{1/k}(A^t_{\tilde{g}}) = f(x).
\end{align}
\end{theorem}
Note that we take the elementary symmetric function 
with respect to the metric $\tilde{g}$.
As noted previously in \cite{Jeff2}, the 
$C^2$ estimate does not work  for $t = 1$, which 
is why we must make the restriction $t<1$.
Also, for $t>1$, the equation is not necessarily elliptic, 
therefore $t=1$ is critical for several reasons. It 
is an interesting problem to investigate the 
limiting behaviour of the solutions in
Theorem \ref{sharpthm} as $t \rightarrow 1$.

 We next write this curvature equation as 
a partial differential equation. We have the following 
formula for the transformation of $A^t$
under a conformal change of metric $\tilde{g} = e^{2w} g$:
\begin{align}
\label{change1}
\tilde{A}^t &= A^t -\nabla^2 w -\frac{1-t}{n-2}(\Delta w)g
+ dw \otimes dw - \frac{2-t}{2} |\nabla w|^2 g.
\end{align}
Since $A^t = A^1 + \frac{1-t}{n-2} tr(A^1) g$,
this formula follows easily from 
the standard formula for the transformation of
the Schouten tensor (see \cite{Jeff2}):
\begin{align}
\tilde{A}^1 = A^1 - \nabla^2 w + dw \otimes dw
- (1/2)| \nabla w|^2 g.
\end{align}
From (\ref{change1}), we may write (\ref{eqn1}) with 
respect to the background metric $g$
\begin{align}
\label{PDE}
 \sigma_k^{1/k} \left( \nabla^2 w + \frac{1-t}{n-2}(\Delta w)g
+ \frac{2-t}{2} |\nabla w|^2 g - dw \otimes dw 
- A^t_g \right) = - f(x) e^{2w} > 0. 
\end{align}

 In the following sections we will derive a priori 
estimates for solutions of (\ref{PDE}), culminating in the 
existence proof in Section \ref{existence}.
In Section \ref{positivecase}, we make some remarks on the 
positive curvature case. 

  We also point out that the Hessian equation
\begin{align}
\sigma_k^{1/k}( \nabla^2 w + S ) = \Psi(x,w)
\end{align}
where $S$ is a symmetric tensor was considered 
in \cite{Yanyan1}, \cite{Delanoe}. For a general existence 
theorem, one must assume that the background metric has non-negative 
sectional curvature. We emphasize that, because of special 
properties of (\ref{PDE}), we do not need such an assumption. 
\subsection{Acknowledgements}
A substantial part of the research for this paper was carried
out while both authors were visiting the Forschungsinstitut f\"ur Mathematik at ETH,
Zurich.  The authors wish to acknowledge the generous support of the Institute, and
the hospitality of Professor Michael Struwe. The authors would 
also like to thank Pengfei Guan for helpful suggestions. 
\section{Ellipticity}
\label{ellipticity}
In this section we will discuss the ellipticity 
properties of equation (\ref{PDE}). 
\begin{definition}
\label{Newtontensor} Let $A : V \rightarrow V$ be 
a symmetric linear transformation,
where $V$ is an $n$-dimensional inner product space.
For $0 \leq q \leq n$, the {\em{$q$th Newton transformation}}
associated with $A$ is
\begin{align*}
T_q(A) = \sigma_q(A) \cdot I - \sigma_{q-1}(A)\cdot A + \cdots +
(-1)^qA^q.
\end{align*}
\end{definition}
Several useful identities involving the Newton tranformation and
elementary symmetric polynomials 
were derived in  
\cite{Reilly}.  First, if $A^i_j$ are the components
of $A$ with respect to some basis of $V$, then
\begin{align}
\label{localNewtontensor} T_q(A)^i_j = \frac{1}{q!} \delta^{i_1
\dots i_q i}_{j_1 \dots j_q j} A_{i_1}^{j_1} \cdots A_{i_q}^{j_q}.
\end{align}
Here $\delta^{i_1 \dots i_q i}_{j_1 \dots j_q j}$ is the
generalized Kronecker delta symbol, and we are
using the Einstein summation convention. Also,
\begin{align}
\label{localA}
\sigma_k(A) = \frac{1}{k!} \delta^{i_1 \dots i_k}_{j_1 \dots j_k}
A_{i_1}^{j_1} \cdots A_{i_k}^{j_k}.
\end{align}
We note that if $A : \mathbf{R} \rightarrow
\mbox{Hom}(V,V)$, then
\begin{align}
\label{dsigmak}
\frac{d}{dt} \sigma_k(A(t)) = T_{k-1}(A(t))^i_j \frac{d}{dt}
A(t)^j_i = T_{k-1} (A(t))^{ij} \frac{d}{dt} A(t)_{ij},
\end{align}
that is, the ($k-1$)-Newton transformation is what we get when we
differentiate $\sigma_k$.
We also note the identities
\begin{align}
\label{TA}
\begin{split}
T_{k-1}(A)^i_j A^j_i &= k \sigma_k(A), \\
tr T_{k-1}(A) &= T_{k-1}(A)^k_k = (n - k + 1)\sigma_k(A).
\end{split}
\end{align}
The following Proposition describes
some important properties of the cones $\Gamma_k^{\pm}$, and their
relation to the Newton transformations:

\begin{proposition}
\label{ellsumm}
(i)  Each set $\Gamma_k^{+}$ is an open convex cone
with vertex at the origin, and we
have the following sequence of inclusions
$$\Gamma_n^{+} \subset \Gamma_{n-1}^{+} \subset \dots \subset
\Gamma_{1}^{+}.$$ \\
(ii) If the eigenvalues of $A$ are in $\Gamma_k^{+}$ (resp., $\Gamma_k^{-}$), then 
$T_{k-1}$ is positive (resp., negative) definite. \\     
(iii)  For symmetric linear transformations $A \in \Gamma_k^+$,
$B \in \Gamma_k^+$, and $t \in [0,1]$,
we have the following inequality
\begin{align}
\label{convexity}
\{ \sigma_k( (1-t)A + tB) \}^{1/k}
\geq (1-t) \{ \sigma_k(A) \}^{1/k}
+ t \{ \sigma_k(B) \}^{1/k}.
\end{align}
\end{proposition}
\begin{proof}
The proof of this proposition is standard, and may be found in
\cite{CNSIII} and \cite{Garding}. 
\end{proof}
Note that the inequality (\ref{convexity})
states that $\sigma_k^{1/k}$ is a {\em{concave}} function 
in $\Gamma_k^+$. This will be essential 
in proving the $C^2$ and $C^{2, \alpha}$ estimates 
in later sections. 
\begin{proposition}
\label{incone}
If $w$ is a solution of (\ref{PDE}) with $A^t_g \in \Gamma_k^-$
for some $t \leq 1$, then 
\begin{align*}
\bar{\nabla}^2 w &\in \Gamma_k^+,  \mbox{  where}\\
\bar{\nabla}^2 w &\equiv  \nabla^2 w + \frac{1-t}{n-2}(\Delta w)g
+ \frac{2-t}{2}|\nabla w|^2 g - dw \otimes dw 
- A^t_g 
\end{align*}
\end{proposition}
\begin{proof}
Since $M$ is compact, at a minimum of the solution $w$
we have
\begin{align*}
\sigma_k^{1/k} \left( \nabla^2 w (p) + \frac{1-t}{n-2} \Delta w(p) g - A^t_g(p)
 \right) = -f(p)e^{2w(p)}> 0,
\end{align*} 
with $\nabla^2w(p)$ positive semidefinite, and therefore 
$ \nabla^2 w (p) + \frac{1-t}{n-2} \Delta w(p) g \geq 0$.
From Proposition \ref{ellsumm}, we conclude 
that $\bar{\nabla}^2 w(p) \in \Gamma_k^{+}$.
Since the cones are connected,
by continuity we have $\bar{\nabla}^2 w \in \Gamma_k^+$.
\end{proof}
\begin{proposition}
\label{ellprop}
If $A_g^t \in \Gamma_k^-$ for some $t \leq 1$, 
then equation {\em{(\ref{PDE})}} is elliptic at any solution.
\end{proposition}
\begin{proof}
We define 
\begin{align*}
F_t[w, \nabla w, \nabla^2 w]
=  \sigma_k \Big( \nabla^2 w  &+ \frac{1-t}{n-2}(\Delta w)g
+ \frac{2-t}{2}|\nabla w|^2 g - dw \otimes dw - A^t_g \Big) \\
&-  |f(x)|^k e^{2kw},
\end{align*}
so that solutions of (\ref{PDE}) are exactly the zeroes of $F_t$. 
We then suppose that $w \in C^2(M) $ satisfies 
$F_t[w, \nabla w, \nabla^2 w] = 0$.  Define $w_s = w + s \varphi$,
then 
\begin{align}
\label{linear}
\begin{split}
\mathcal{L}^t(\varphi) &= \frac{d}{ds}F_t[w_s, \nabla w_s, \nabla^2 w_s] \Big|_{s=0}\\
&= \frac{d}{ds}\left( \sigma_k ( \bar{\nabla}^2 w_s) \right) \Big|_{s=0}
- \frac{d}{ds} \left( |f(x)|^k e^{2kw_s} \right) \Big|_{s=0}.
\end{split}
\end{align}
From (\ref{dsigmak}), we have
\begin{align*}
\frac{d}{ds} \left( \sigma_k (  \bar{\nabla}^2 w_s) \right) \Big|_{s=0}
= T_{k-1} (  \bar{\nabla}^2 w )_{ij}  (\bar{\nabla}^2 w_s)'_{ij}.
\end{align*}
where
\begin{align*}
(\bar{\nabla}^2 w_s)' = \frac{d}{ds} \left( \bar{\nabla}^2 w_s \right)\Big|_{s=0}.
\end{align*}
We compute
\begin{align*}
(\bar{\nabla}^2 w_s)' = \nabla^2 \varphi  + \frac{1-t}{n-2}(\Delta \varphi)g
+ (2-t)\langle dw, d \varphi \rangle g - 
2 dw \otimes d \varphi.
\end{align*}
Therefore,
\begin{align}
\notag
\frac{d}{ds} \left( \sigma_k (  \bar{\nabla}^2 w_s) \right) \Big|_{s=0}
= T_{k-1} (  \bar{\nabla}^2 w)_{ij} 
\Big(  \nabla_i \nabla_j \varphi  & + \frac{1-t}{n-2}(\Delta \varphi)g_{ij}\\
\label{1stterm}
&+ (2-t) \langle dw, d \varphi \rangle g_{ij} - 
2 w_i \varphi_j \Big).
\end{align}
For the second term on the $RHS$ of (\ref{linear}) we have
\begin{align}
\label{2ndterm}
\frac{d}{ds} \left( |f(x)|^k e^{2kw_s} \right) \Big|_{s=0}
= 2k|f(x)|^k e^{2k w} \varphi.
\end{align}
Combining (\ref{1stterm}) and (\ref{2ndterm}), we conclude
\begin{align}
\label{ker}
\mathcal{L}^t(\varphi) = T_{k-1}( \bar{\nabla}^2w)_{ij}( \nabla_i \nabla_j \varphi 
+ \frac{1-t}{n-2} \Delta \varphi g_{ij})  -  2k|f(x)|^k e^{2k w} \varphi
+ \cdots
\end{align}
where $+ \cdots$ denotes additional terms which are linear in $\nabla \varphi$.
Defining
\begin{align}
\label{Qdef}
(Q_{k-1})_{ij}= (T_{k-1})_{ij} + \frac{1-t}{n-2} (T_{k-1})_{ll} \delta_{ij},
\end{align}
we have 
\begin{align}
\label{ker2}
\mathcal{L}^t(\varphi) =  Q_{k-1}( \bar{\nabla}^2w)_{ij}( \nabla_i
\nabla_j \varphi)  -  2k|f(x)|^k e^{2k w} \varphi
+ \cdots
\end{align}
By Propositions \ref{ellsumm} and  \ref{incone}, 
$Q_{k-1}( \bar{\nabla}^2w)_{ij}  > 0$, so $\mathcal{L}^t$ is elliptic.  
\end{proof}
\noindent
Since the coefficient of $\varphi$ in the 
zeroth-order term of (\ref{ker2}) is strictly 
negative, we have 
\begin{corollary}
\label{linvert}
If $A_g^t \in \Gamma_k^-$ for some $t \leq 1$, then at any 
solution of {\em{(\ref{PDE})}}, the linearized operator 
$\mathcal{L}^t : C^{2, \alpha}(M) \rightarrow C^{\alpha}(M)$
is invertible.
\end{corollary}
\section{$C^0$ estimate}
\label{C0}
 We begin with an important property of $\sigma_k$ in the cone $\Gamma_k^{+}$.
\begin{lemma}
\label{sks}
Let $A$ and $B$ be symmetric $n \times n$ matrices. 
Assume that $A$ is positive semi-definite, $B \in \Gamma_k^+$, 
and $A + B \in \Gamma_k^+$. Then 
$$ \sigma_k(A + B) \geq \sigma_k (B).$$
If $A$ is negative semi-definite, then 
$$ \sigma_k(A + B) \leq \sigma_k (B).$$
\end{lemma}
\begin{proof}
Let $F(t) = \sigma_k (tA + B) - \sigma_k(B)$ for $t \in [0,1]$.
Note that from convexity of the cone $\Gamma_k^+$, 
we have $ t(A+ B) + (1-t) B = tA + B \in \Gamma_k^+$. 
Using (\ref{dsigmak}), we have 
\begin{align*}
F'(t) = T_{k-1}(tA + B)^{ij} A_{ij}  \geq 0,
\end{align*}
since $T_{k-1}(tA + B)$ is positive definite from 
Proposition \ref{ellsumm}.
Therefore $F(t)$ is non-decreasing, and $F(0) = 0$, so we have
$F(1) = \sigma_k(A + B) - \sigma_k(B) \geq 0.$
The negative case is similar. 
\end{proof}
\begin{proposition} 
\label{easyC0}
Suppose $A^t_g \in \Gamma_k^-$ for some $t \leq 1$. Then there exist 
constants $  \underline{\delta} < 0 < \overline{\delta}$
depending only upon $f$, $A_g^t$ and $k$, 
such that for any solution $w(x)$ of {\em{(\ref{PDE})}}, we have 
$\underline{\delta} < w(x) < \overline{\delta}$.
\end{proposition}
\begin{proof}
Since $N$ is compact, at a minimum of the function $w(x)$
we have
\begin{align*}
\sigma_k^{1/k} \left( \nabla^2 w (p) + \frac{1-t}{n-2} \Delta wg(p) - A_g^t(p)
 \right) = - f(p) e^{w(p)}
\end{align*} 
with $\nabla^2 w (p)$ positive semidefinite, and therefore 
$ \nabla^2 w (p) + \frac{1-t}{n-2} \Delta w g(p) \geq 0$. 
From Proposition \ref{incone} and Lemma \ref{sks}  we have 
\begin{align*}
\sigma_k^{1/k} (-A_g^t(p)) \leq - f(p) e^{w(p)},
\end{align*}
and certainly we can choose $\underline{\delta}$ such that
\begin{align*}
w(x) \geq w(p) \geq \mbox{ln} \left( \frac{ \sigma_k^{1/k}(-A_g^t(p))}{-f(p)} \right)
\geq \mbox{ln} \left( \underset{x \in N}{ \mbox{min}} 
\frac{ \sigma_k^{1/k}(-A_g^t(x))}{-f(x)} \right) > \underline{\delta}.
\end{align*}
Similarly, if the maximum of $w(x)$ is at $q \in N$, we
can choose  $\overline{\delta}$ such that
\begin{align*}
w(x) \leq w(q) \leq \mbox{ln} \left( \frac{ \sigma_k^{1/k}(-A_g^t(q))}{-f(q)} \right)
\leq \mbox{ln} \left( \underset{x \in N}{ \mbox{max}} 
\frac{ \sigma_k^{1/k}(-A_g^t(x))}{-f(x)} \right) < \overline{\delta} .
\end{align*}
\end{proof}
\section{$C^1$ estimate}
\label{C1}
\begin{proposition}
\label{C1estimate}
 Let $w$ be a $C^3$ solution of (\ref{PDE}) for some $t \leq 1$, 
satisfying $\underline{\delta} < w < \overline{\delta}$.
Then $\Vert \nabla w \Vert_{L^{\infty}} < C_2$, where $C_2$ 
depends only upon $\underline{\delta},\overline{\delta}, g, t$. 
\end{proposition}
We consider the following function
\begin{align*}
h = \left( 1 + \frac{|\nabla w|^2}{2} \right) e^{\phi(w)},
\end{align*}
where $\phi : \mathbf{R} \rightarrow \mathbf{R}$ is a 
function of the form 
\begin{align*}
\phi(s) = c_1(c_2 + s)^p.
\end{align*}
The constants $c_1, c_2,$ and $p$ will be chosen later. 
We will estimate the maximum value of the 
function $h$, and this will give us the gradient 
estimate. 

  Since $N$ is compact, and $h$ is 
continuous, we suppose the maximum of $h$ occurs and a point 
$p \in N$. We take a normal coordinate system 
$(x^1, \dots, x^n)$ at $p$. Then we have 
$g_{ij}(p) = \delta_{ij}$, and $\Gamma^i_{jk}(p) = 0$, 
where $g = g_{ij} dx^i dx^j$, and $\Gamma^i_{jk}$
is the Christoffel symbol (see \cite{Besse}). 

 Locally, we may write $h$ as 
\begin{align*}
h = \left( 1 + \frac{1}{2} g^{lm}w_{l}w_{m} \right)
e^{\phi(w)} = v e^{\phi(w)}.
\end{align*}
In a neighborhood of $p$, differentiating $h$ in the $x^i$
direction we have
\begin{align}
\notag
\partial_i h &= h_i = \frac{1}{2} \partial_i ( g^{lm}w_l w_m)
e^{\phi(w)} + v e^{\phi(w)}\phi'(w)w_i\\
\label{goog}
& = \frac{1}{2} \partial_i ( g^{lm})w_l w_m e^{\phi(w)} 
 + g^{lm} \partial_i(w_l) w_m e^{\phi(w)} + v e^{\phi(w)}\phi'(w)w_i 
\end{align}
Since in a normal coordinate system, the first 
derivatives of the metric vanish at $p$, 
and since $p$ is a maximum for $h$, 
evaluating (\ref{goog}) at $p$, we have
\begin{align}
\label{hip}
w_{li}w_l = - v \phi'(w) w_i.
\end{align}
Next we differentiate (\ref{goog}) in the $x^j$ direction. 
Since $p$ is a maximum, $\partial_j \partial_i h = h_{ij}$ is negative 
semidefinite, and we get (at $p$)
\begin{align*}
0 \gg h_{ij} &= \frac{1}{2} \partial_j \partial_i g^{lm} w_l w_m
e^{\phi(w)} + w_{lij}w_l e^{\phi(w)} + w_{li}w_{lj} e^{\phi(w)}
+ w_{li}w_{l} e^{\phi(w)} \phi'(w) w_j \\
&+ v_j e^{\phi(w)} \phi'(w) w_i
 + v e^{\phi(w)} (\phi'(w))^2 w_i w_j + v e^{\phi(w)} \phi''(w)w_j w_i
+ v e^{\phi(w)} \phi'(w) w_{ij}
\end{align*} 
Next we note that $v_j = w_{lj} w_l$, and 
using (\ref{hip}), we have 
\begin{align*}
0 \gg h_{ij} &= \frac{1}{2} \partial_j \partial_i g^{lm} w_l w_m
e^{\phi(w)} + w_{lij}w_l e^{\phi(w)} + w_{li}w_{lj} e^{\phi(w)}\\
&+ ( \phi''(w) - \phi'(w)^2) v e^{\phi(w)} w_i w_j 
 + v e^{\phi(w)} \phi'(w) w_{ij}.
\end{align*} 
We recall from Section \ref{ellipticity} that
\begin{align}
\label{Qagain}
\bar{Q}^t_{ij} = \bar{T}_{ij} + \frac{1-t}{n-2} \bar{T}_{ll}
\delta_{ij}, 
\end{align}
is positive definite, where $\bar{T}_{ij}$
means $T_{k-1}(\bar{\nabla}^2 w)_{ij}$.

So we divide by $v e^{\phi(w)}$, sum 
with $\bar{Q}^t_{ij}$, and we have
the inequality
\begin{align}
\label{inequality1}
0 \geq \frac{1}{2v} \bar{Q}^t_{ij} \partial_i \partial_j
g^{lm}w_l w_m + \frac{1}{v} \bar{Q}^t_{ij} w_{lij}w_l
+ ( \phi''(w) - \phi'(w)^2)\bar{Q}^t_{ij} w_i w_j
+ \phi'(w) \bar{Q}^t_{ij} w_{ij},
\end{align}
since $w_{li} w_{lj}$ is positive semidefinite.

 We will use equation (\ref{PDE}) to replace the $w_{ij}$ term
with lower order terms,  and then differentiate equation
(\ref{PDE}) in order to replace the $w_{lij}$ term with 
lower order terms. With respect to our local coordinate system, 
from (\ref{change1}) we have 
\begin{align}
\notag
\bar{\nabla}^2w_{ij} = &- A^t_{ij} + w_{ij} - w_r \Gamma^r_{ij} 
+\frac{1-t}{n-2}(w_{kk} - w_r \Gamma^r_{kk})g_{ij}
- w_i w_j \\
\label{wbarlocal}
& + \frac{2-t}{2} (g^{r_1 r_2} w_{r_1} w_{r_2} )g_{ij}.
\end{align} 
At the point $p$, this simplifies to 
\begin{align}
\label{wijeqn}
\bar{\nabla}^2w_{ij} = - A^t_{ij}
+ w_{ij}  +\frac{1-t}{n-2}(w_{kk})g_{ij}
- w_i w_j + \frac{2-t}{2} (|\nabla w|^2 )\delta_{ij}.
\end{align}
We write equation (\ref{PDE}) in our local coordinate system
\begin{align}
\label{localnegeqn}
 \sigma_k( g^{lj} (\bar{\nabla}^2 w)_{ij} ) = 
|f(x)| e^{2kw}.
\end{align}
Note that the $g^{lj}$ term is present since we need 
to raise an index on the tensor before we apply $\sigma_k$.
From (\ref{wijeqn}), we have at $p$,
\begin{align*}
\bar{Q}^t_{ij}w_{ij} &= \bar{T}_{ij}\left(  \bar{\nabla}^2w_{ij}
+ A^t_{ij} -\frac{1-t}{n-2}(\Delta w )g_{ij}
+ w_i w_j - \frac{2-t}{2} (| \nabla w|^2 )\delta_{ij}\right)\\
 & \ \ \ \ \ \ \ 
+ \frac{1-t}{n-2}\Delta w \bar{T}_{ll}\\
& = \bar{T}_{ij}\left( \bar{\nabla}^2w_{ij}
+ A^t_{ij}
+ w_i w_j - \frac{2-t}{2} (| \nabla w|^2 )\delta_{ij}\right).
\end{align*}
From the identities in (\ref{TA}), we have
\begin{align*}
\bar{T}_{ij}(\bar{\nabla}^2w)_{ij}
= k \sigma_k( \bar{\nabla}^2w).
\end{align*}
Therefore, using equation (\ref{localnegeqn}), we have 
\begin{align}
\notag
\bar{Q}^t_{ij}w_{ij} 
& =  k \sigma_k( \bar{\nabla}^2w)
+ \bar{T}_{ij}\left( A^t_{ij}
+ w_i w_j - \frac{2-t}{2} (| \nabla w|^2 )\delta_{ij}\right)\\
\label{Qijform}
& = k |f(x)| e^{2kw}+ \bar{T}_{ij}\left( A^t_{ij}
+ w_i w_j - \frac{2-t}{2} (| \nabla w|^2 )\delta_{ij}\right).
\end{align}
Next we take $m$ with $1 \leq m \leq n$, and apply $\partial_m$
to (\ref{localnegeqn})
\begin{align}
\label{firstpartial}
\partial_m \left\{ \sigma_k ( g^{lj} \bar{\nabla}^2w_{ij}  ) \right\} = 
\partial_m (|f(x)| e^{2kw}).
\end{align}
For the left side we have
\begin{align*}
& \partial_m \left\{  \sigma_k  ( g^{lj} \bar{\nabla}^2w_{ij} ) \right\} = 
\bar{T}_{il} \partial_m ( g^{lj}\bar{\nabla}^2w_{ij} )
=\bar{T}_{il} \left( 
\partial_m g^{lj}\bar{\nabla}^2w_{ij} 
+  g^{lj} \partial_m \bar{\nabla}^2w_{ij} \right)\\
&= \bar{T}_{il} 
\Bigg( 
\partial_m g^{lj} \bar{\nabla}^2w_{ij}+  g^{lj} \partial_m 
\Bigg\{ - A^t_{ij}
+ w_{ij} - w_r \Gamma^r_{ij} 
+\frac{1-t}{n-2}(w_{kk} - w_r \Gamma^r_{kk})g_{ij}\\
&\ \ \ \ \ \ \ \ \ - w_i w_j + \frac{2-t}{2} 
(g^{r_1 r_2} w_{r_1} w_{r_2} )g_{ij}  \Bigg\} \Bigg)\\
&= \bar{T}_{il} \Bigg( 
\partial_m g^{lj}  \bar{\nabla}^2w_{ij}+  
g^{lj} \Bigg\{ - \partial_m A^t_{ij}
+ w_{ijm} - w_{rm} \Gamma^r_{ij} - w_r \partial_m  \Gamma^r_{ij} \\
&+\frac{1-t}{n-2}(w_{kkm} - w_{rm} \Gamma^r_{kk} - w_{r} \partial_m \Gamma^r_{kk})g_{ij} 
+\frac{1-t}{n-2}(w_{kk} - w_r \Gamma^r_{kk}) \partial_m g_{ij}
- 2w_{im} w_j \\
&+ \frac{2-t}{2} \Big( (\partial_m g^{r_1 r_2}) w_{r_1} w_{r_2} g_{ij}
+ 2 g^{r_1 r_2} w_{r_1m} w_{r_2} g_{ij}+ 
g^{r_1 r_2} w_{r_1} w_{r_2} \partial_m g_{ij}  \Big) \Bigg\} \Bigg).
\end{align*}
Evaluating this expression at $p$, we have
\begin{align}
\notag
\partial_m \left\{  \sigma_k  ( g^{lj} \bar{\nabla}^2w_{ij} ) \right\}
&= \bar{T}_{ij} \Bigg( 
-\partial_m A^t_{ij}
+ w_{ijm} - w_r \partial_m  \Gamma^r_{ij} 
+\frac{1-t}{n-2}(w_{kkm} - w_{r} \partial_m \Gamma^r_{kk})g_{ij} \\
\label{fdp}
& \ \ \ \ \ \ \ \ \ - 2w_{im} w_j + (2-t) ( 
 w_{lm} w_{l}) g_{ij}  \Bigg).
\end{align}
For the right hand side of (\ref{firstpartial}) we have
\begin{align}
\label{rhs1d}
\partial_m \left( |f(x)| e^{2kw}\right)
&= \partial_m |f(x)| e^{2kw} 
+ 2k |f(x)| e^{2kw} w_m.
\end{align}
We have obtained the expansion of (\ref{firstpartial}):
\begin{align}
\begin{split}
\label{fullfirstder}
&\bar{T}_{ij} \Bigg( 
-\partial_m A^t_{ij}
+ w_{ijm} - w_r \partial_m  \Gamma^r_{ij} 
+\frac{1-t}{n-2}(w_{kkm} - w_{r} \partial_m \Gamma^r_{kk})\delta_{ij} \\
& \ \ - 2w_{im} w_j - (2-t) v \phi'(w) w_m  \delta_{ij} 
\Bigg)
= (\partial_m |f(x)|) e^{2kw}
+ 2k |f(x)| e^{2kw}w_m. 
\end{split}
\end{align}
Note that the third order terms in the above expression are 
\begin{align*}
&\bar{T}_{ij} \left( w_{ijm} +\frac{1-t}{n-2}w_{kkm}\delta_{ij} \right)
= \bar{Q}^t_{ij} w_{ijm}.    
\end{align*}
Next we sum (\ref{fullfirstder}) with $w_m$, using (\ref{hip})
we have the following formula
\begin{align}
\begin{split}
\label{firstderterm}
& \bar{Q}_{ij} w_m w_{ijm} + 
\bar{T}_{ij} \Bigg( 
 -w_m \partial_m A^t_{ij} - w_m w_r \partial_m  \Gamma^r_{ij} 
-\frac{1-t}{n-2}( w_{r} w_m \partial_m \Gamma^r_{kk})\delta_{ij} \\
& \ \ \ \ \ \ \ \ \ + 2 v \phi'(w) w_i w_j - (2-t) 
v \phi'(w) | \nabla w|^2  \delta_{ij} \Bigg)\\
& \hspace{30mm} = w_m (\partial_m f(x)) e^{2kw} 
+ 2k |f(x)| e^{2kw} | \nabla w|^2. 
\end{split}
\end{align}
Substituting (\ref{Qijform}) and (\ref{firstderterm}) into 
(\ref{inequality1}), we arrive at the inequality
\begin{align*}
0 &\geq \frac{1}{2v} \bar {Q}^t_{ij} \partial_i \partial_j
 g^{lm}w_l w_m 
 + ( \phi''(w) - \phi'(w)^2)\bar{Q}^t_{ij} w_i w_j\\
&+ \frac{1}{v} 
\bar{T}_{ij} \Bigg( 
w_m \partial_m A^t_{ij} + w_m w_r \partial_m  \Gamma^r_{ij} 
+\frac{1-t}{n-2}( w_{r} w_m \partial_m \Gamma^r_{kk})\delta_{ij}
 - 2 v \phi'(w) w_i w_j\\
& + (2-t) v \phi'(w) | \nabla w|^2  
\delta_{ij} \Bigg)
+\frac{1}{v} \Big( w_m (\partial_m |f(x)|) e^{2kw} 
+ 2k |f(x)| e^{2kw} | \nabla w|^2  \Big)\\
&+ \phi'(w) \left( k |f(x)| e^{2kw}
+  \bar{T}_{ij} 
\Big( 
A^t_{ij} + w_i w_j - 
\frac{2-t}{2} (| \nabla w|^2 )\delta_{ij}
\Big) \right).
\end{align*}
Recalling that 
$\bar{Q}^t_{ij} = \bar{T}_{ij} + \frac{1-t}{n-2} \bar{T}_{ll}\delta_{ij}$, 
we have
\begin{align}
\notag
0 &\geq 
 \Big( \phi''(w) - \phi'(w)^2 - \phi'(w) \Big) 
\bar{T}_{ij}w_i w_j\\
\notag
&+ \frac{1}{v} 
\bar{T}_{ij} \Bigg( \frac{1-t}{n-2}  \partial_k \partial_k
g^{lm}w_l w_m (\delta_{ij}/2)
+ \frac{1-t}{n-2} v( \phi''(w) - \phi'(w)^2)| \nabla w|^2 g_{ij}\\
\notag
&+ \frac{1}{2} \partial_i \partial_j g^{lm}w_l w_m 
+ w_m \partial_m A^t_{ij} + w_m w_r \partial_m  \Gamma^r_{ij} 
+\frac{1-t}{n-2}( w_{r} w_m \partial_m \Gamma^r_{kk})\delta_{ij}\\
\notag
&+ v \phi'(w) A^t_{ij} 
+ \frac{2-t}{2} v \phi'(w)| \nabla w|^2 \delta_{ij} \Bigg)\\
\label{mess}
&+ \frac{1}{v} \Big( w_m (\partial_m f(x)) e^{2kw}
+ 2k |f(x)| e^{2kw} | \nabla w|^2
\Big) + \phi'(w) \left( 
  k|f(x)| e^{2kw}\right).
\end{align}
\begin{lemma}
At $p$, in normal coordinates, we have
\begin{align*}
 \sum_{l,m} (\partial_i \partial_j
g^{lm} + 2 \partial_l \Gamma^m_{ij}) u_l u_m 
= 2 \sum_{l,m} R_{iljm}u_l u_m,
\end{align*}
where $R_{iljm}$ are the components of the Riemann
curvature tensor of $g$.
\end{lemma}
\begin{proof}
For the proof, see \cite{Jeff2}. 
\end{proof}
Using the lemma, and collecting terms in (\ref{mess}), we arrive at
\begin{align}
\notag
0 &\geq 
 \Big( \phi''(w) - \phi'(w)^2 - \phi'(w) \Big) 
\bar{T}_{ij}w_i w_j\\
\notag
&+ \bar{T}_{ij} \Bigg( \frac{1-t}{n-2}  R_{klkm} \frac{w_l w_m}{v} \delta_{ij}
+ R_{iljm} \frac{w_l w_m}{v}
+ \frac{1-t}{n-2} ( \phi''(w) - \phi'(w)^2)|\nabla w|^2 g_{ij}\\
\notag
& \hspace{30mm} + \frac{w_m}{v} \partial_m A^t_{ij} + \phi'(w) A^t_{ij} 
+ \frac{2-t}{2} \phi'(w)| \nabla w|^2 \delta_{ij} \Bigg)\\
\label{mess2}
& \ \ \ \ \ \ + \frac{1}{v} \Big( w_m (\partial_m f(x)) e^{2kw}
+ 2k |f(x)| e^{2kw} | \nabla w|^2
\Big) + \phi'(w) \left( 
  k|f(x)| e^{2kw}\right).
\end{align}
Now we will choose $\phi(s)$. 
\begin{lemma}
\label{choosephi}
Assume that $ \underline{\delta} < s < \overline{\delta}$.
Then we may choose constants $c_1, c_2,$ and $p$
depending only upon $\underline{\delta}$, and $\overline{\delta}$.
so that $\phi(s) = c_1(c_2 +s)^p$ satisfies
\begin{align}
\label{phi'}
 \phi'(s) > 0,
\end{align}
and 
\begin{align}
\label{phi''}
 \phi''(s) - \phi'(s)^2 - \phi'(s) > 0.
\end{align}
\end{lemma}
\begin{proof}
We have 
\begin{align*}
\phi'(s) = pc_1(c_2 +s)^{p-1},
\end{align*}
and
\begin{align*}
\phi''(s) = p(p-1) c_1 (c_2 +s)^{p-2}.
\end{align*}
To satisfy (\ref{phi'}) we need $c_1 > 0$, $p>0$, and
$c_2 + s> 0$. So choose $c_2 >  \underline{\delta}$.
Next we have
\begin{align*}
\phi''(s) - \phi'(s)^2 - \phi'(s) &= 
 p(p-1) c_1 (c_2 +s)^{p-2} - ( pc_1(c_2 +s)^{p-1})^2
- p c_1(c_2 +s)^{p-1}\\
&= pc_1 (c_2 + s)^{p-2} \Big( (p-1) - pc_1 (c_2 +s)^p
- (c_2 +s) \Big).
\end{align*}
Now choose 
\begin{align*}
c_1 = \frac{1}{p^2 \cdot \mbox{ max} \{ (c_2 +s)^p \} },
\end{align*}
and $p$ so large that 
\begin{align*}
 \underline{\delta} < c_2 < - \overline{\delta} + p -1 - \frac{1}{p}.
\end{align*}
Then we have
\begin{align*}
\phi''(s) - \phi'(s)^2 - \phi'(s) &\geq 
 \frac{1}{ p \cdot \mbox{ max} \{ (c_2 +s)^p \} } (c_2 +s)^{p-2}
\Big( p - 1 - \frac{1}{p} - c_2 - s \Big)\\
& > \frac{1}{ p \cdot \mbox{ max} \{ (c_2 +s)^p \} } (c_2 +s)^{p-2}
( \overline{\delta} - s ) >0.
\end{align*} 
\end{proof}
With $\phi(s)$ chosen as above, for  $s \in [ \underline{\delta}, 
\overline{\delta}]$, we let 
$$\epsilon_1 = \mbox{min}\{ \phi'(s) \},$$ and
$$\epsilon_2 = \mbox{min} \{ \phi''(s) - \phi'(s)^2 - \phi'(s) \}.$$  
From the inequality (\ref{mess2}), we have
\begin{align}
\notag
0 &\geq 
 \epsilon_2 \bar{T}_{ij}w_i w_j 
+ \bar{T}_{ij} \Bigg( \frac{1-t}{n-2}  R_{klkm} \frac{w_l w_m}{v} \delta_{ij}
 + R_{iljm} \frac{w_l w_m}{v}\\
\notag
& \hspace{30mm} + \frac{w_m}{v} \partial_m A^t_{ij} + \phi'(w) A^t_{ij} 
+ \Big( \frac{2-t}{2} + \frac{1-t}{n-2} \Big) 
\epsilon_1 | \nabla w|^2 \delta_{ij} \Bigg)\\
\label{mess3}
& \ \ \ \ \ \ + \frac{1}{v} \Big( w_m (\partial_m f(x)) e^{2kw}
+ 2k |f(x)| e^{2kw} | \nabla w|^2
\Big) + \epsilon_1 \left( 
  k|f(x)| e^{2kw}\right).
\end{align}
If the matrix 
\begin{align}
\begin{split}
\label{zash}
&\Big( \frac{1-t}{n-2}  R_{klkm} \delta_{ij}
+ R_{iljm} \Big) \frac{w_l w_m}{v}
+ \frac{w_m}{v} \partial_m A^t_{ij} + \phi'(w) A^t_{ij} 
+ \Big( \frac{2-t}{2} + \frac{1-t}{n-2} \Big) 
\epsilon_1 | \nabla w|^2 \delta_{ij} 
\end{split}
\end{align}
has an eigenvalue less than $1$, then the gradient estimate 
is immediate since the last term dominates
($\frac{2-t}{2} + \frac{1 -t}{n-2} > 0$). 
Otherwise, absorbing lower order terms in (\ref{mess3}), we have 
\begin{align}
\label{thingy2}
C &\geq 
 \epsilon_2 \bar{T}_{ij}w_i w_j + \bar{T}_{ll}.
\end{align}
From (\ref{TA}) we obtain 
\begin{align*}
  \sigma_{k-1} \leq C.
\end{align*}
\begin{proposition}
\label{compactness}
Let $k \geq 2$, and $A \in \Gamma_k^+$ be a symmetric 
linear transformation. If $0 < c_1 \leq \sigma_k(A)$,
and $\sigma_{k-1}(A) \leq c_2$, then we have a bound on 
the eigenvalues of $A$, that is, $|\lambda(A)| \leq C$,
where $C$ depends only on $c_1$ and $c_2$.
\end{proposition}
\begin{proof}
The proof may be found in \cite{Yanyan1}.
\end{proof}
Using this result, if $k \geq 2$, we see that 
\begin{align*}
|\lambda| \leq C,
\end{align*}
and since $\bar{T}_{k-1}$ is positive definite, this
implies (see \cite{Yanyan1})
\begin{align*}
\bar{T}_{k-1}^{ii} \geq \frac{1}{C}>0, \mbox{ for } i = 1 \dots n.  
\end{align*}
Equation (\ref{thingy2}) then implies that
\begin{align*}
|\nabla w|^2 \leq C.
\end{align*}
Note that in the case $k=1$, we do not require the proposition since
$T_0^{ij} = \delta^{ij}$, and therefore 
(\ref{thingy2}) gives the gradient bound.
\section{$C^2$ estimate}
\label{C2}
\begin{proposition}
\label{C2estimate}
 Let $w$ be a $C^4$ solution of (\ref{PDE}) for some $t < 1$ 
satisfying $\underline{\delta} < w < \overline{\delta}$,
and $\Vert \nabla w \Vert_{L^{\infty}} < C_1$.
Then  $\Vert \nabla^2 w \Vert_{L^{\infty}} \leq C_2$,
where $C_2$ 
depends only upon $\underline{\delta},\overline{\delta}, 
 C_1, g, t$. 
\end{proposition}
Let $S(TN)$ denote the unit tangent bundle of $N$, and 
we consider the following function $h: S(TN) \mapsto \mathbf{R}$,
\begin{align*}
h(e_p) = ( \nabla^2 w + \Lambda |\nabla w|^2 g)(e_p, e_p),
\end{align*}
where $\Lambda$ is a constant to be chosen later. 
Since $S(TN)$ is compact, 
let $h$ have a maximum at the vector $\tilde{e}_p$.
We use normal coordinates at $p$, and by rotating,
assume that the tensor is diagonal at $p$, 
and without loss of generality, we may assume 
that $\tilde{e}_p = {\partial}/{\partial x^1}$, and 
that $\nabla^2 w$ is diagonal at $p$. 
  
 We let $\tilde{h}$ denote the function defined 
in a neighborhood of $p$
\begin{align*}
\tilde{h}(x) &= ( \nabla^2 w + \Lambda | \nabla w|^2 g)( {\partial}/{\partial
 x^1}, {\partial}/{\partial x^1} )\\
& = (\nabla^2 w)_{11} + \Lambda |\nabla w|^2 \\
& = w_{11} - \Gamma^l_{11}w_l + \Lambda |\nabla w|^2.
\end{align*}
Differentiating in the $i$th coordinate direction, we obtain
\begin{align}
\label{wix}
\tilde{h}_i = w_{11i} - \partial_i \Gamma^l_{11} w_l
- \Gamma^l_{11} w_{li} + \Lambda \partial_i g^{kl}w_k w_l
+ 2\Lambda g^{kl} (\partial_i w_k) w_l.
\end{align}
The function $\tilde{h}(x)$ has a maximum at $p$,
so evaluating (\ref{wix}) at $p$, we obtain
\begin{align}
\label{wip}
w_{11i} =  \partial_i \Gamma^l_{11} u_l + 2 \Lambda w_{ik}w_{k}.
\end{align}
Next we differentiate (\ref{wix}) in the $x^j$ direction. 
Since $p$ is a maximum, $\partial_j \partial_i \tilde{h} = \tilde{h}_{ij}$
is negative semidefinite, and we get (at $p$)
\begin{align*}
0 \gg \tilde{w}_{ij} & = w_{11ij} - \partial_i \partial_j
\Gamma^l_{11} w_l - \partial_i \Gamma^l_{11} w_{lj}
- \partial_j \Gamma^l_{11} u_{li}\\
& + \Lambda \partial_j \partial_i g^{kl}w_k w_l 
+ 2 \Lambda w_{ikj}w_k + 2 \Lambda w_{ik}w_{kj}.
\end{align*}
We again recall from Section \ref{ellipticity} that
\begin{align}
\label{Qagain2}
\bar{Q}^t_{ij} = \bar{T}_{ij} + \frac{1-t}{n-2} \bar{T}_{ll}
\delta_{ij}, 
\end{align}
is positive definite, where $\bar{T}_{ij}$
means $T_{k-1}(\bar{\nabla}^2 w)_{ij}$.
We sum with $\bar{Q}^t_{ij}$ and we have the inequality
\newcommand{\bq}{\bar{Q}^t_{ij}}
\begin{align}
\label{testineq}
\begin{split}
0 \geq & \bq w_{11ij} - \bq \partial_i \partial_j
\Gamma^l_{11} w_l - 2\bq  \partial_i \Gamma^l_{11} w_{lj}
+ 2 \Lambda \bq  \partial_j \partial_i g^{kl}w_k w_l  \\
& + 2 \Lambda \bq  w_{ik}w_{kj} + 2 \Lambda \bq  w_{ijk}w_k 
\end{split}
\end{align}
We will use (\ref{fdp}) to replace the last term, and 
we will differentiate equation (\ref{PDE}) twice 
to replace the first term. 

 Using (\ref{fdp}) and (\ref{rhs1d}), we have at $p$, 
\begin{align}
\begin{split}
\label{fdp2}
\bar{T}_{ij} &\Big( 
-\partial_m A^t_{ij}- w_r \partial_m  \Gamma^r_{ij} 
-\frac{1-t}{n-2}( w_{r} \partial_m \Gamma^r_{kk})g_{ij} 
- 2w_{im} w_j + (2-t) ( 
 w_{lm} w_{l}) g_{ij}  \Big)\\
& \ \ \ \ \ \ \ \ \ \ + \bar{Q}^t_{ij}w_{ijm}   = \partial_m |f(x)| e^{2kw} 
+ 2k |f(x)| e^{2kw} w_m.
\end{split}
\end{align}
  The next step is to rewrite the second derivative 
terms in terms of $ \bar{\nabla}^2 w$. 
To further simplify notation, we let $\bar{w}_{ij} = (\bar{\nabla}^2 w)_{ij}$.
We have
\begin{align}
\label{nice1}
 w_{ij} = \bar{w}_{ij} - \frac{1-t}{n-2} \Delta w \delta_{ij}
 - \frac{2-t}{2} |\nabla w|^2 \delta_{ij}
+ w_i w_j + (A^t_g)_{ij}.
\end{align}
Since we are in the cone $\Gamma_k^+$, the trace is
positive by Proposition \ref{ellsumm},
and since $\bar{w}_{11}$ is the largest eigenvalue,
we have
\begin{align}
\label{nice2}
 |\bar{w}_{ii}| \leq (n-1) \bar{w}_{11}, \ \ \ \ i = 1 \dots n.
\end{align}
Using (\ref{nice1}) and (\ref{nice2}), 
we may estimate (\ref{fdp2}) 
\begin{align}
\label{fdp3}
\bar{Q}^t_{ij}w_{ijm}w_m &
\geq C + C \sum_i \bar{T}_{ii} + C \bar{w}_{11} \sum_i \bar{T}_{ii}.
\end{align}
Inequality (\ref{testineq}) may then be rewritten as
\begin{align}
\label{testineq2}
\begin{split}
0 \geq & \bq w_{11ij} 
+ C + C \sum_i \bar{T}_{ii} + C \bar{w}_{11} \sum_i \bar{T}_{ii}.
+ \frac{2 \Lambda(1-t)}{n-2} \bar{w}_{11}^2 \sum_i \bar{T}_{ii}. 
\end{split}
\end{align}
We recall that the equation is 
$$\sigma_k^{1/k}( \bar{\nabla}^2 w) = -f(x) e^{2w}.$$
To simplify notation, write $ \sigma = \sigma_k^{1/k}$.
Differentiating once in the $x^1$ direction, we have
\begin{align*}
\frac{\partial \sigma}{\partial r_{ij}}( \partial_1 (\bar{\nabla}^2 w)^i_j)
= - f_1 e^{2w} - 2f e^{2w}w_{1}.
\end{align*}
Differentiating twice, we obtain
\begin{align*}
&\partial_1 \left(\frac{\partial \sigma }{\partial r_{ij}} \right)
( \partial_1 (\bar{\nabla}^2 w)^i_j)
+ \frac{\partial \sigma }{\partial r_{ij}}
(\partial_1 \partial_1 (\bar{\nabla}^2 w)^i_j)\\
&= \Big( \frac{\partial^2 \sigma}{\partial r_{ij} \partial r_{lm}} \Big)
( \partial_1 (\bar{\nabla}^2 w)^l_m)
( \partial_1 (\bar{\nabla}^2 w)^i_j)
+ \frac{\partial \sigma }{\partial r_{ij}}
(\partial_1 \partial_1 (\bar{\nabla}^2 w)^i_j)\\
&= - f_{11} e^{2w} -4f_1 e^{2w} w_i 
- 4 f e^{2w} w_1^2 - 2f e^{2w} w_{11}.
\end{align*}
Since $\sigma_k^{1/k}$ is concave in $\Gamma_k^+$,
we have the inequality
\begin{align}
\label{ineq1}
\bar{T}_{k-1}^{ij} \Big( \partial_1 \partial_1 (\bar{\nabla}^2 w)^i_j \Big)
\geq (-f(x)e^{2w})^{k-1} ( - f_{11} e^{2w} -4f_1 e^{2w} w_i 
- 4 f e^{2w} w_1^2 - 2f e^{2w} w_{11}).
\end{align}
Using our assumptions, we have 
\begin{align}
\label{easyineq1}
\bar{T}_{k-1}^{ij} \Big( \partial_1 \partial_1 (\bar{\nabla}^2 w)^i_j \Big)
\geq C + C \bar{w}_{11}.
\end{align}
We differentiate (\ref{wbarlocal}) twice, and evaluate 
at $p$ to obtain 
\begin{align*}
\partial_1 \partial_1 & (\bar{\nabla}^2 w)_{ij} = - A^t_{ij11} + w_{ij11} 
+ 2 w_{r1} \partial_1(\Gamma^r_{ij})  + w_r \partial_1 
 \partial_1(\Gamma^r_{ij}) \\
& + \frac{1-t}{n-2} (w_{kk11})g_{ij} + \frac{1-t}{n-2} w_{kk} \partial_1 
 \partial_1 (g_{ij}) \\
& - \frac{1-t}{n-2} ( 2 w_{r1} \partial_1 \Gamma^r_{kk} + w_r 
 \partial_1 \partial_1 \Gamma^r_{kk}) 
- ( w_{i11}w_j + 2 w_{i1} w_{j1} + w_i w_{j11})  \\
& + \frac{2-t}{2} \Big( \partial_1 \partial_1 
g^{r_1 r_2} w_{r_1} w_{r_2} \delta_{ij}
 + 2 (w_{r11} w_r + w_{r1}w_{r1}) \delta_{ij}
+ |\nabla w|^2   \partial_1 \partial_1 (g_{ij}) \Big).
\end{align*} 
From (\ref{wip}) we can replace terms of the form 
$w_{11i}$ and we have
\begin{align}
\label{bigmess}
\begin{split}
&\bar{T}_{k-1}^{ij} \bar{w}_{ij11} 
= \bq \bar{w}_{ij11} + 
\bar{T}_{k-1}^{ij} \Big\{ - A^t_{ij11} + 2 w_{r1} 
\partial_1(\Gamma^r_{ij})  + w_r \partial_1 
 \partial_1(\Gamma^r_{ij}) \\
& + \frac{1-t}{n-2} w_{kk} \partial_1 
 \partial_1 (g_{ij}) 
- \frac{1-t}{n-2} ( 2 w_{r1} \partial_1 \Gamma^r_{kk} + w_r 
 \partial_1 \partial_1 \Gamma^r_{kk})
+ 4 \Lambda w_{ik} w_k w_j - 2 w_{i1} w_{j1}   \\
& + \frac{2-t}{2} \Big( \partial_1 \partial_1 
g^{r_1 r_2} w_{r_1} w_{r_2} \delta_{ij}
 -4 \Lambda w_{rm}w_m w_r \delta_{ij} 
+2 w_{r1}w_{r1}) \delta_{ij}
+ |\nabla w|^2   \partial_1 \partial_1 (g_{ij}) \Big) \Big\}
\end{split}
\end{align}
Substituting (\ref{bigmess}) in (\ref{easyineq1}), we have
\begin{align}
\label{bigmess2}
\begin{split}
\bq w_{ij11} \geq 
 C + C \bar{w}_{11} + C \sum_i \bar{T}_{ii} 
+ C \bar{w}_{11} \sum_i \bar{T}_{ii}
- (2-t) \bar{w}_{11}^2 \sum_i \bar{T}_{ii}. 
\end{split}
\end{align}
Next we substitute inequality (\ref{bigmess2})
into (\ref{testineq2}) and we obtain
\begin{align}
\begin{split}
\label{est}
0 \geq & C + C \bar{w}_{11} 
+ C \sum_i \bar{T}_{ii} + C \bar{w}_{11} \sum_i \bar{T}_{ii}
+ \left( \frac{ 2 \Lambda (1-t)}{n-2} - (2-t) \right) \bar{w}_{11}^2 \sum_i 
\bar{T}_{ii}. 
\end{split}
\end{align}
Since $t<1$, we may choose $\Lambda$ large to dominate the
$-(2-t)$ term (this is the point where the assumption 
$t<1$ is crucial). Choosing 
\begin{align*}
\Lambda \geq \frac{n-2}{1-t}\left( 1 + \frac{2-t}{2} \right),
\end{align*}
we obtain
\begin{align}
\label{coolineq}
 C + C \bar{w}_{11} + C \sum_i \bar{T}_{k-1}^{ii} 
+ C \bar{w}_{11} \sum_i \bar{T}_{k-1}^{ii} \geq  
2 \bar{w}_{11}^2 \sum_i \bar{T}_{k-1}^{ii}.
\end{align}
Dividing by $2 \bar{w}_{11}^2$ and using (\ref{TA}), we obtain
\begin{align}
\label{secder}
\begin{split}
\sigma_{k-1}  \leq 
\Big( \frac{C_1}{ \bar{w}_{11}^2} + \frac{C_2}{\bar{w}_{11}}
\Big) \sigma_{k-1} + \frac{C}{ \bar{w}_{11}^2} 
+ \frac{C}{\bar{w}_{11}}. 
\end{split}
\end{align}
If  
\begin{align*}
\frac{C_1}{ \bar{w}_{11}^2} + \frac{C_2}{\bar{w}_{11}}  \geq \frac{1}{2},
\end{align*}
then we have the necessary eigenvalue bound.
So we may assume that
\begin{align*}
\frac{C_1}{ \bar{w}_{11}^2} + \frac{C_2}{\bar{w}_{11}}
  \leq \frac{1}{2},
\end{align*}
and substitution into inequality (\ref{secder}) yields
\begin{align*}
\frac{1}{2} \sigma_{k-1}  \leq 
\frac{C}{ \bar{w}_{11}^2} 
+ \frac{C}{\bar{w}_{11}}. 
\end{align*}
Without loss of generality we may assume that $\bar{w}_{11} \geq 1$,
and from the above inequality we obtain
\begin{align*}
\sigma_{k-1} \leq C,
\end{align*}
which by Proposition \ref{compactness} yields the 
eigenvalue bound in the case $k \geq 2$. 
In the case $k=1$, (\ref{secder}) already gives the 
eigenvalue estimate. 
\section{Proof of Theorem \ref{sharpthm}}
\label{existence}

  In the previous sections we have demonstrated an a priori 
$C^2$ estimate for solutions of (\ref{PDE}). Therefore 
(\ref{PDE}) is uniformly elliptic with respect to 
any solution. It is straightforward to 
verify that (\ref{PDE}) is a concave function of the 
second derivative variables. It follows from the work  
of Evans \cite{Evans}, and Krylov \cite{Krylov} 
that there exists a constant $C$ such that
for any solution $w$ of (\ref{PDE}) we have 
$$ \parallel \! w  \! \! \parallel_{C^{2,\alpha}} < C.$$
For $s \in [0,1]$, we consider the equation
\begin{align}
\begin{split}
\label{PDE2}
\sigma_k^{1/k} & \left( \nabla^2 w^s + \frac{1-t}{n-2}(\Delta w^s)g
+ \frac{2-t}{2}|\nabla w^s|^2 g - dw^s \otimes dw^s - H(s) \right) \\
& =   \Big( s|f(x)| + (1-s) \Big) e^{2w^s} > 0, 
\end{split}
\end{align}
where $H(s) = s A^t_g - (1-s) { {n}\choose{k}}^{-1/k}g$. 
Define
$$S = \{ s \in [0,1] : (\ref{PDE2}) \mbox{ has a solution } 
w^s \in C^{2, \alpha} \}.$$ 
It follows from Corollary \ref{linvert} and the implicit 
function theorem (see \cite{GT}) that $S$ is open. 
For $s = 0$, we have the solution 
$w \equiv 0$, therefore $S$ is non-empty. 
Clearly, the estimates from the previous sections remain valid
upon replacing $A^t_g$ by $H(s)$, 
and $|f(x)|$ by $ s|f(x)| + (1-s) $.  
Consequently there 
exists a constant $C$ independent of $s$ such that
$$ \parallel \! w^s  \! \! \parallel_{C^{2,\alpha}} < C,$$
which implies that $S$ is closed. From connectedness, $S = [0,1]$,
therefore there exists a solution at $s=1$.

  Using the maximum principle as in Proposition \ref{easyC0},
it follows that $w \equiv 0$ is the unique solution at
$s = 0$. To prove uniqueness at $s=1$, assume by contradiction 
that we have solutions $w_1$ and $w_2$ for $s = 1$. 
Then we may run the continuity method in reverse, 
starting at $s = 1$ and descending to $s=0$. From the 
a priori estimates, we obtain 2 paths of solutions,
one starting at $w_1$, and another starting at $w_2$. 
Since there is a unique solution at $s=0$, the paths must 
coincide at some $s \geq 0$. This contradicts 
local invertibility. 
 
\section{Remarks on the positive curvature case}
\label{positivecase}
In this paper, we have concentrated on the case 
that $A^t_g \in \Gamma_k^-$, but it is also 
interesting to consider the case 
$A^t_g \in \Gamma_k^+$. This problem was studied 
in \cite{Jeff2} for $t=1$, where 
the compactness was reduced to an $L^{\infty}$
estimate. The estimate proved in this paper 
also reduce the compactness to an $L^{\infty}$ 
estimate for $ - \infty < t \leq 1$. 

For $f(x) > 0$, and background metric $g$ with 
$A^t_g \in \Gamma_k^+$, we consider a conformal change of 
metric $\tilde{g} = e^{-2w} g$. Then (\ref{eqn1}) becomes
\begin{align}
\label{PDE3}
 \sigma_k^{1/k} \left( \nabla^2 w + \frac{1-t}{n-2}(\Delta w)g
- \frac{2-t}{2}|\nabla w|^2 g + dw \otimes dw 
+ A^t_g \right) = f(x) e^{2w} > 0. 
\end{align}
Note the gradient terms have a different sign now,
but the estimates for higher derivatives still work. 
In particular, for the $C^1$ estimate, the only 
modification necessary is the following variant 
of Lemma \ref{choosephi}
\begin{lemma}
\label{choosephipos}
Assume that $ \underline{\delta} < s < \overline{\delta}$.
Then we may choose constants $c_1, c_2,$ and $p$
depending only upon $\underline{\delta}$, and $\overline{\delta}$.
so that $\phi(s) = c_1(c_2 +s)^p$ satisfies
\begin{align}
\label{phi'2}
 \phi'(s) < 0,
\end{align}
and 
\begin{align}
\label{phi''2}
 \phi''(s) - \phi'(s)^2 - \phi'(s) > 0.
\end{align}
\end{lemma}
The proof given in Section \ref{C1} then works as before
for all $t \leq 1$. 
For the $C^2$ estimate, using the following test function 
\begin{align*}
h(e_p) = ( \nabla^2 w + \Lambda |\nabla w|^2 g + dw \otimes dw)(e_p, e_p),
\end{align*}
the estimate (\ref{est}) becomes
\begin{align}
\begin{split}
\label{est2}
0 \geq C + C \bar{w}_{11} + 
C \sum_i \bar{T}_{ii} + C \bar{w}_{11} \sum_i \bar{T}_{ii}
+ \left( \frac{2 \Lambda (1-t)}{n-2} + (2-t) \right) \bar{w}_{11}^2 \sum_i 
{T}_{ii}. 
\end{split}
\end{align}
The estimate now works for all $ - \infty < t \leq 1$. 
Note we can include the endpoint $t=1$ in the 
positive case. 

We conclude with an outline of some progress that has been 
made for $t=1$ in the positive case. 
For $\sigma_2$ in dimension 4, the $L^{\infty}$ 
estimate has been proved (if $M$ is not conformally 
equivalent to $S^4$) in  
Chang, Gursky and Yang (see \cite{CGYnew}).
Existence of solutions in the 
locally conformally flat case has 
been demonstrated in \cite{LiLi} and \cite{GuanWang}.  
A sufficient condition for the $L^{\infty}$ estimate 
in the determinant case may be found in \cite{Jeff2}. 


\bibliography{references}
\noindent
\small{\textsc{Department of Mathematics, University of Notre Dame, 
Notre Dame, IN 46556}}\\
{\em{E-mail Address:}} \ {\texttt{mgursky@nd.edu}}\\
\\
\small{\textsc{Department of Mathematics, Massachusetts Institute
of Technology, Cambridge, MA 02139}}\\
{\em{E-mail Address:}} \
{\texttt{jeffv@math.mit.edu}}
\end{document}